\documentclass[12pt,letterpaper]{article}
\usepackage{ifthen,latexsym,amssymb,amsmath}

\newcommand{\C}[1]{{\protect\cal #1}}
\newcommand{\B}[1]{{\bf #1}}
\newcommand{\I}[1]{{\mathbb #1}}
\renewcommand{\O}[1]{\overline{#1}}
\newcommand{\ceil}[1]{\lceil #1\rceil}
\newcommand{\e}{\varepsilon}
\newcommand{\floor}[1]{\lfloor #1\rfloor}
\newcommand{\me}{{\mathrm e}}
\renewcommand{\mid}{:}

\newif\ifnotesw\noteswtrue
\newcommand{\comment}[1]{\ifnotesw $\blacktriangleright$\ {\sf #1}\ 
  $\blacktriangleleft$ \fi}
\noteswfalse	

\newcommand{\beq}[1]{\begin{equation}\label{eq:#1}}
\newcommand{\eeq}{\end{equation}}
\newcommand{\req}[1]{\textrm{(\ref{eq:#1})}}

\newtheorem{theorem}{Theorem}
\newcommand{\bth}[2][nothing]{\ifthenelse{\equal{#1}{nothing}}
 {\begin{theorem}} {\begin{theorem}[#1]}\label{th:#2}}

\newtheorem{lemma}[theorem]{Lemma}
\newcommand{\blm}[2][nothing]{\ifthenelse{\equal{#1}{nothing}}
 {\begin{lemma}} {\begin{lemma}[#1]}\label{lm:#2}}

\newtheorem{problem}[theorem]{Problem}
\newcommand{\bpr}[2][nothing]{\ifthenelse{\equal{#1}{nothing}}
 {\begin{problem}} {\begin{problem}[#1]}\label{pr:#2}}

\newcommand{\case}[1]{\smallskip\noindent{\bf Case #1} }
\newcommand{\claim}[1]{\smallskip\noindent{\bf Claim #1} }

\newcommand{\bpf}[1][Proof.]{\smallskip\par\noindent{\it #1} }
\newcommand{\qed}{\nolinebreak\mbox{\hspace{5 true pt}%
  \rule[-0.85 true pt]{3.9 true pt}{8.1 true pt}}}
\newcommand{\cqed}{\nolinebreak\mbox{\hspace{5 true pt}%
  \rule[-0.85 true pt]{2.0 true pt}{8.1 true pt}}}
\newcommand{\epf}{\qed \medskip}
\newcommand{\bcpf}{\bpf[Proof of Claim.]}
\newcommand{\ecpf}{\cqed \medskip}

\newcommand{\brm}{\smallskip\noindent{\bf Remark.} }

\begin{document}

\newcommand{\brac}[1]{\left(#1\right)}
\newcommand{\bfrac}[2]{\brac{\frac{#1}{#2}}}
\newcommand{\set}[1]{\{#1\}}
\newcommand{\stack}[2]{\genfrac{}{}{0pt}{}{#1}{#2}}
\def\cT{{\cal T}}
\def\bM{{\bf M}}
\def\bk{{\bf k}}
\def\cC{{\cal C}}
\def\cE{{\cal E}}
\def\cF{{\cal F}}
\def\qwert{2^{\tilde O(1/\ep^2)}}
\def\cond{{\rm Cond}}
\def\hu{\hat{u}}
\def\bC{{\bf C}}
\def\bT{{\bf T}}
\def\hDD{\hat{{\bf D}}}
\def\bM{{\bf M}}
\def\DD{{\bf D}}
\def\norm{\hbox{Norm}}
\def\half{{\textstyle{1\over2}}}
\def\bR{{\bf R}}
\def\recip#1{{1\over#1}}
\def\dij{\hat{d}_{i,j}}
\def\vol{{\rm Vol}}
\def\for{\hbox{  for  }}
\def\hT{\hat{T}}
\def\bA{{\bf A}}
\def\bW{{\bf W}}
\def\bhW{\hat{\bW}}
\def\hc{\hat{c}}
\def\ep{\epsilon}
\def\ind{{\rm ind}}
\def\hw{\hat{w}}
\def\bs{\bar{\sigma}}
\def\A{{\cal A}}
\def\B{{\cal B}}
\def\F{{\cal F}}
\def\E{{\bf E}}
\def\H{{\cal H}}
\def\Ex{{\bf E}}
\def\cS{{\cal S}}
\def\cP{{\cal P}}
\def\cQ{{\cal Q}}
\def\ind{{\rm ind}}   
\def\bn{\bar{n}}
\def\bv{\bar{\n}}
\def\a{\alpha}
\def\b{\beta}
\def\d{\delta}
\def\D{\Delta}
\def\e{\epsilon}
\def\f{\phi}
\def\g{\gamma}
\def\G{\Gamma}
\def\k{\kappa}
\def\la{\lambda}
\def\K{\Kappa}
\def\z{\zeta}
\def\th{\theta}
\def\TH{\Theta}
\def\l{\lambda}
\def\m{\mu}
\def\n{\nu}
\def\p{\pi}
\def\P{\Pi}
\def\r{\rho}
\def\R{\Rho}
\def\s{\sigma}
\def\S{\Sigma}
\def\t{\tau}
\def\om{\omega}
\def\OM{\Omega}
\def\bigmid{\rule[-3.5mm]{0.1mm}{9mm}}
\def\N{{\cal N}}
\def\Pr{\mbox{{\bf Pr}}}
\def\CG{{\cal G}}
\def\CA{{\cal A}}
\def\whps{{\bf whp }}
\def\whp{{\bf whp}}
\def\Prob{{\bf Pr}}
\def\he{\hat{\e}}
\def\BD{{\bf D}}
\def\bW{{\bf W}}
\def\bB{{\bf B}}
\def\hr{\hat{r}}
\def\hR{\hat{R}}

\newcommand{\ratio}[2]{\mbox{${#1\over #2}$}}
\newcommand{\bbD}[1]{\bar{{\bf D}}^{(#1)}}
\newcommand{\gap}[1]{\mbox{\hspace{#1 in}}}
\newcommand{\bD}[1]{{\bf D}^{(#1)}}
\newcommand{\hbD}[1]{\hat{{\bf D}}^{(#1)}}
\newcommand{\bTT}[1]{{\bf T}^{(#1)}}
\newcommand{\limninf}{\mbox{$\lim_{n \rightarrow \infty}$}}
\newcommand{\proofstart}{{\bf Proof\hspace{2em}}}
\newcommand{\tset}{\mbox{$\cal T$}}
\newcommand{\proofend}{\hspace*{\fill}\mbox{$\Box$}}
\newcommand{\bfm}[1]{\mbox{\boldmath $#1$}}
\newcommand{\reals}{\mbox{\bfm{R}}}
\newcommand{\expect}{\mbox{\bf E}}
\newcommand{\Exp}{\mbox{\bf E}}
\newcommand{\card}[1]{\mbox{$|#1|$}}
\newcommand{\scaps}[1]{\mbox{\sc #1}}
\newcommand{\rdup}[1]{\lceil #1 \rceil }
\newcommand{\rdown}[1]{\lfloor #1 \rfloor }
\newcommand{\mnote}[1]{\marginpar{\footnotesize\raggedright#1}}
\newcommand{\rt}{\right}
\newcommand{\lt}{\left}

\newcommand{\re}{\mbox{\rm e}}
\newcommand{\setm}{\setminus}

\newenvironment{proof}{\noindent{\bf Proof\,}}{\hfill$\Box$}
\newtheorem{remark}{Remark}

\def\bx{{\bf x}}
\def\gNM2{{\cal G}_{\n,\m}^{\d\geq 2}}
\def\GNM2{{G}_{\n,\m}^{\d\geq 2}}



\newcommand{\con}[2]{#1\leftrightarrow #2}
\newcommand{\dist}{\mathrm{dist}}
\newcommand{\DT}{D^{\mathrm{tree}}}
\newcommand{\ehr}{\mbox{\sc Ehr}}
\newcommand{\diam}{\mathrm{diam}}

\title{First Order Definability of Trees and\\ Sparse Random Graphs}

\renewcommand{\thefootnote}{\fnsymbol{footnote}}

\author{Tom Bohman\footnotemark[1]\ \footnotemark[2], Alan Frieze\footnotemark[1]\ \footnotemark[3],
  Tomasz {\L}uczak\footnotemark[4], Oleg
Pikhurko\footnotemark[1]\ \footnotemark[5],\\ Clifford
Smyth\footnotemark[1], Joel
 Spencer\footnotemark[6], and Oleg Verbitsky\footnotemark[7]}

\date{}

\maketitle

\footnotetext[1]{Department of Mathematical Sciences,
Carnegie Mellon University,
Pittsburgh, PA 15213, USA.}
\footnotetext[2]{Partially supported by NSF grant DMS-0401147.}
\footnotetext[3]{Partially supported by NSF Grant CCR-0200945.}
\footnotetext[4]{Department of Discrete Mathematics, Adam Mickiewicz
  University, Pozna\'n 61-614, Poland. Partially supported by KBN grant 
1 P03A 025 27.}
\footnotetext[5]{Partially supported by the Berkman Faculty
Development Fund, CMU.}
\footnotetext[6]{Courant Institute, New York University, New York, NY
10012, USA.}
\footnotetext[7]{Institut f\"ur Informatik, Humboldt Universit\"at,
Berlin 10099, Germany.
Supported by an Alexander von Humboldt fellowship.}

\renewcommand{\thefootnote}{\arabic{footnote}}

\vspace{5mm}

\begin{abstract}
 Let $D(G)$ be the smallest quantifier depth of a first order formula
which is true for a graph $G$ but false for any other non-isomorphic
graph. This can be viewed as a measure for the first order descriptive
complexity of $G$.

We will show that almost surely $D(G)=\Theta(\frac{\ln n}{\ln\ln n})$,
where $G$ is a random tree of order $n$ or the giant component of a
random graph $\C G(n,\frac cn)$ with constant $c>1$. These results
rely on computing the maximum of $D(T)$ for a tree $T$ of order $n$
and maximum degree $l$, so we study this problem as
well.\end{abstract}

\section{Introduction}

This paper deals with graph properties expressible in first order
logic. The vocabulary consists of variables, connectives ($\vee$,
$\wedge$ and $\neg$), quantifiers ($\exists$ and $\forall$), and two
binary relations: the equality and the graph adjacency ($=$ and $\sim$
respectively). The variables denote vertices only so we are not
allowed to quantify over sets or relations. The notation $G\models A$
means that a graph $G$ is a model for a \emph{sentence} $A$ (a first
order formula without free variables); in other words, $A$ is true
for the graph $G$. All sentences and graphs are assumed to be
finite. The Reader is referred to Spencer's book~\cite{spencer:slrg}
(or to~\cite{kim+pikhurko+spencer+verbitsky:03rsa}) for more details.

A first order sentence $A$ \emph{distinguishes} $G$ from $H$ if
$G\models A$ but $H\not\models A$. Further, we say that $A$
\emph{defines} $G$ if $A$ distinguishes $G$ from any non-isomorphic
graph $H$. In other words, $G$ is the unique (up to an isomorphism)
finite model for $A$.

The \emph{quantifier depth} (or simply \emph{depth}) $D(A)$ is the
largest number of nested quantifiers in $A$. This parameter is closely
related to the complexity of checking whether $G\models A$. 

The main parameter we will study is $D(G)$, the smallest quantifier
depth of a first order formula defining $G$. It was first
systematically studied by Pikhurko, Veith and
Verbitsky~\cite{pikhurko+veith+verbitsky:03} (see
also~\cite{pikhurko+verbitsky:03}). In a sense, a defining formula $A$
can be viewed as the canonical form for $G$ (except that $A$ is not
unique): in order to check whether $G\cong H$ it suffices to check
whether $H\models A$. Unfortunately, this approach does not seem to
lead to better isomorphism algorithms but this notion, being on the
borderline of combinatorics, logic and computer science, is
interesting on its own and might find unforeseen applications.

Within a short time-span various results on the values of $D(G)$ for
order-$n$ graphs appeared. The initial
papers~\cite{pikhurko+veith+verbitsky:03,pikhurko+verbitsky:03}
studied the maximum of $D(G)$ (the `worst' case).  The `best' case is
considered by Pikhurko, Spencer, and
Verbitsky~\cite{pikhurko+spencer+verbitsky:04} while Kim, Pikhurko,
Spencer and Verbitsky~\cite{kim+pikhurko+spencer+verbitsky:03rsa}
obtained various results for random graphs.

Here we study these questions for trees and sparse random
structures. Namely, the three main questions we consider are:
 \begin{description}
 \item[Section~\ref{general}:] What is $\DT(n,l)$, the maximum of
$D(T)$ over all trees of order at most $n$ and maximum degree at most
$l$?
 \item[Section~\ref{giant}:] What is $D(G)$, where $G$ is the giant
component of a random graph $\C G(n,\frac{c}{n})$ for constant
$c>1$?
 \item[Section~\ref{random}:] What is $D(T)$ for a random tree $T$ of
order $n$?
 \end{description}

In all cases we determine the order of magnitude of the studied
function.  Namely, we prove that $\DT(n,l)=\Theta(\frac{l\ln n}{\ln
l})$, and whp we have $D(G)=\Theta(\frac{\ln n}{\ln\ln n})$, whenever
$G$ is a random tree of order $n$ or the giant component of a random
graph $\C G(n,\frac cn)$ with constant $c>1$. (The acronym \emph{whp}
stands for `with high probability', i.e.,\ with probability $1-o(1)$.)
Moreover, for some cases involving trees we estimate the smallest
quantifier depth of a first order formula defining $G$ up to a factor
of $1+o(1)$. For instance, we show that for a random tree $T$ of order
$n$ we have whp $D(T)=(1+o(1))\frac{\ln n}{\ln\ln n}$.
 \comment{
also we prove that $\DT(n,l)=(1/2+o(1))\frac{l\ln n}{\ln l}$
whenever both $l=l(n)$ and $\ln n/\ln l$ tends to
infinity as $n\to\infty$.
 }%

\section{Further Notation and Terminology}

Our main tool in the study of $D(G)$ is the \emph{Ehrenfeucht
game}. Its description can be found in Spencer's book~\cite{spencer:slrg}
whose terminology we follow (or
see~\cite[Section~2]{kim+pikhurko+spencer+verbitsky:03rsa}), so here
we will be very brief.

Given two graphs $G$ and $G'$, the \emph{Ehrenfeucht game}
$\ehr_k(G,G')$ is a perfect information game played by two players,
called \emph{Spoiler} and \emph{Duplicator}, and consists of $k$
rounds, where $k$ is known in advance to both players. For
brevity, let us refer to Spoiler as `him' and to Duplicator as
`her'. In the $i$-th round, $i=1,\dots,k$, Spoiler selects one of the
graphs $G$ and $G'$ and marks one of its vertices by $i$; Duplicator
must put the same label $i$ on a vertex in the other graph. At the end
of the game let $x_1,\dots,x_k$ be the vertices of $G$ marked
$1,\dots,k$ respectively, regardless of who put the label there; let
$x_1',\dots,x_k'$ be the corresponding vertices in $G'$. Duplicator
wins if the correspondence $x_i\leftrightarrow x_i'$ is a partial
isomorphism, that is, we require that $\{x_i,x_j\}\in E(G)$ iff
$\{x_i',x_j'\}\in E(G')$ as well as that $x_i=x_j$ iff
$x_i'=x_j'$. Otherwise, Spoiler wins.

The key relation is that $D(G,G')$, the smallest depth of a first
order sentence $A$ distinguishing $G$ from $G'$, is equal to the
smallest $k$ such that Spoiler can win $\ehr_k(G,G')$. Also,
 \beq{D}
 D(G)=\max_{G'\not\cong G} D(G,G'),
 \eeq
 see e.g.~\cite[Lemma~1]{kim+pikhurko+spencer+verbitsky:03rsa}.

Sometimes it will be notationally more convenient to prove the bounds on
$D(G,G')$ for colored graphs which generalize the usual (uncolored)
graphs.  Graphs $G,G'$ are \emph{colored} if we have unary relations
$U_i:V(G)\cup V(G')\to\{0,1\}$, $i\in I$. We say that the vertices in
the set $U_i^{-1}(1)$ have color $i$. Note that some vertices may be
uncolored and some may have more than one color. There are no
restrictions on a color class, i.e.,\ it does not have to be an
independent set.  When the Ehrenfeucht game is played on colored
graphs, Duplicator must additionally preserve the colors of vertices.

Colorings can be useful even if we prove results for uncolored
graphs. For example, if $x\in V(G)$ and $x'\in V(G')$ were selected in
some round, then, without changing the outcome of the remaining game,
we can remove $x$ and $x'$ from $G$ and $G'$ respectively, provided we
color their neighbors with a new color. (Note that in an optimal
strategy of Spoiler, there is no need to select the same vertex
twice.)

We will also use the following fact, which can be easily deduced from
the general theory of the Ehrenfeucht game. Let $x,y\in V(G)$ be
distinct vertices. Then the smallest quantifier depth of a first order
formula $\Phi(z)$ with one free variable $z$ such that $G\models
\Phi(x)$ but $G\not\models \Phi(y)$ is equal to the minimum $k$ such
that Spoiler can win the $(k+1)$-round game $\ehr_{k+1}(G,G)$, where
the vertices $x_1=x$ and $x_1'=y$ have been selected in the first
round.

In this paper $\ln$ denotes the natural logarithm, while the
logarithm base $2$ is written as $\log_2$.

\section{General Trees}\label{general}

Let $\DT(n,l)$ be the maximum of $D(T)$ over all colored trees of
order at most $n$ and maximum degree at most $l$. We split the
possible range of $l,n$ into a few cases.

\bth{MaxDeg} Let both $l$ and $\ln n/\ln l$ tend to the
infinity. Then
 \beq{MaxDeg}
 \DT(n,l)= \left(\frac12+o(1)\right)\, \frac{ l\ln n}{\ln
l}.
 \eeq
 In fact, the lower bound can be achieved by uncolored trees.
 \end{theorem}

In order to prove Theorem~\ref{th:MaxDeg} we need some preliminary
results. Let $\dist_G(x,y)$ denote the distance in $G$ between $x,y\in
V(G)$.

\blm{Distance} Suppose $x,y\in V(G)$ at distance $k$ were selected while
their counterparts $x',y'\in V(G')$ are at a strictly larger
distance (possibly infinity). Then Spoiler can win in at most
$\log_2k+1$ extra moves, playing all of the time inside
$G$.\end{lemma}
 \bpf We prove the claim by induction on $k$. Assume $k\ge 2$ and
choose an appropriate $xy$-path $P$. Spoiler selects a vertex $w\in
V(G)$ which is a \emph{middle vertex} of $P$, that is,
$k_1=\dist_P(x,w)$ and $k_2=\dist_P(y,w)$ differ at most by
one. Suppose that Duplicator responds with $w'\in G'$. It is
impossible that $G'-z'$ contains both an $x'w'$-path of length at most
$k_1$ and a $y'w'$-path of length at most $k_2$. If, for example, the
latter does not exist, then we apply induction to $y,w\in G$. The
required bound follows by observing that $k_1,k_2\le \ceil{\frac
k2}$.\epf

The same method gives the following lemma.

\blm{path} Let $G,G'$ be colored graphs. Suppose that $x,y\in V(G)$
and $x',y'\in V(G')$ have been selected such that $G$ contains some
$xy$-path $P$ of length at most $k$ such that some vertex of $P$ has
color $c$ while this is not true with respect to $G'$. Then Spoiler
can win in at most $\log_2 k +1$ moves playing all of the time inside $G$.

The same conclusion holds if all internal vertices of $P$ have colors
from some fixed set $A$ while any $x'y'$-path of length at most $k$
has a color not in $A$.\qed\end{lemma}

\blm{Tree} Let $T$ be a tree of order $n$ and let $T'$ be a graph
which is not a tree. Then $D(T,T')\le
\log_2n+3$.\end{lemma}
 \bpf
 If $T'$ is not connected, Spoiler selects two vertices $x',y'\in T'$
from different components. Then he switches to $G$ and applies
Lemma~\ref{lm:Distance}, winning in at most $\log_2 n+3$ moves in
total.

Otherwise, let $C'\subset T'$ be a cycle of the shortest length
$l$. If $l>2n+1$, then Spoiler picks two vertices $x',y'$ at distance
at least $n$ in $C'$ (or equivalently in $T'$). But the diameter of
$T$ is at most $n-1$, Spoiler switches to $T$ and starts halving the
$xy$-path, making at most $\log_2 n+3$ moves in total, cf.\
Lemma~\ref{lm:Distance}.

If $l\le 2n+1$, then Spoiler selects some three adjacent vertices of
$C'$, say $x',z',y'$ in this order. Now, he applies
Lemma~\ref{lm:path} with respect to $k=l-2$.\epf

\bpf[Proof of Theorem~\ref{th:MaxDeg}.] Let us prove the upper bound
first.

Let $T$ be any tree of order at most $n$ and maximum degree at most
$l$. Let $T'$ be an arbitrary colored graph not isomorphic to $T$. By
Lemma~\ref{lm:Tree} we can assume that $T'$ is a tree.

In fact, we will be proving the upper bound on the version of the
$(T,T')$-game, wherein some distinguished vertex, called the
\emph{root}, is given and all graph isomorphisms must additionally
preserve the root. (This can be achieved by introducing a new color
$U_0$ which is assigned to the root only.) The obtained upper bound,
if increased by $1$, applies to the original function $D(T,T')$
because we can regard $x_1$ and $x_1'$, the first two moves of the
Ehrenfeucht game, as the given roots.

It is easy to show that $T$ contains a vertex $x\in T$ such that any
component of $T-x$ has order at most $\frac n2$. We call such a vertex
a \emph{median} of $T$. Spoiler selects this vertex $x$; let
Duplicator reply with $x'$. We can assume that the degrees of $x$ and
$x'$ are the same: otherwise Spoiler can exhibit this discrepancy in
at most $l+1$ extra moves.
 \comment{Alternating sides at most once.}%

We view the components of $T-x$ and $T'-x'$ as colored rooted graphs
with the neighbors of $x$ and $x'$ being the roots. As $T\not\cong
T'$, some component $C_1$ has different multiplicities $m_1$ and
$m_1'$ in $T-x$ and $T'-x'$. As $d(x)=d(x')$, we have at least two
such components. Assume that for $C_1$ and $C_2$ we have $m_1>m_1'$
and $m_2<m_2'$. By the condition on the maximum degree, $m_1'+m_2\le
l-1$. Hence, $\min(m_1',m_2)\le \frac{l-1}2$. Let us assume, for
example, that $m_1'\le \frac{l-1}2$. Spoiler chooses the roots of any
$m_1'+1$ $C_1$-components of $T-x$. It must be the case that some
vertices $y\in V(T)$ and $y'\in V(T')$ have been selected, so that $y$
lies in a $C_1$-component $F\subset T-x$ while $y'$ lies in a
component $F'\subset T'-x$ not isomorphic to $C_1$.  Let $n_1$ be the
number of vertices in $F$. By the choice of $x$, $n_1\le \frac n2$.

Now, Spoiler restricts his moves to $V(F)\cup V(F')$. If Duplicator
moves outside this set, then Spoiler uses Lemma~\ref{lm:path},
winning in at most $\log_2n+O(1)$ moves. Otherwise Spoiler uses the
recursion applied to $F$.

Let $f(n,l)$ denote the largest number of moves (over all trees $T,T'$
with $v(T)\le n$, $\Delta(T)\le l$, and $T\not\cong T'$) that
Duplicator can survive against the above strategy with the additional
restriction that a situation where Lemma~\ref{lm:path} can be applied
never occurs and we always have that $d(x)=d(x')$. Clearly,
 \beq{DTf}
 \DT(n,l)\le f(n,l) + \log_2n + l +O(1).
 \eeq

As $m_1\le \frac{n-1}{n_1}$, we get the following
recursive bound on $f$.
 \beq{DT}
 \textstyle
 f(n,l)\le \max\Big\{2 + \min(\frac{l-1}2,\frac{n-1}{n_1}) +
 f(n_1,l)\mid 1\le n_1\le \frac n2\Big\}.
 \eeq
 Denoting $n_0=n$ and unfolding~\req{DT} as long as $n_i\ge 1$, say $s$
times, we obtain that $f(n,l)$ is bounded by the maximum of 
 \beq{f}
 2s  + \sum_{i=1}^s
\min\left(\frac{l-1}2,\frac{n_{i-1}}{n_i}\right),
 \eeq
 over all sequences $n_1,\dots,n_s$ such that
 \beq{n}
 1\le n_i \le \frac{n_{i-1}}2,\quad i\in[s].
 \eeq
  Note that the restrictions~\req{n} force $s$ to be at most $\log_2
n$. Let us maximize~\req{f} over all $s\in\I N$
and real $n_i$'s satisfying~\req{n}. 

It is routine to see that for the optimal sequence we have $2\le
\frac{n_{i-1}}{n_i}\le \frac{l-1}2$, $i\in[s]$; moreover, both these
inequalities can be simultaneously strict for at most one index
$i$. 
 \comment{Indeed, suppose on the contrary that for two indexes $1\le
i<j< s$ we have $2<n_i/n_{i+1}<\frac{l-1}2$ and
$2<n_j/n_{j+1}<\frac{l-1}2$. Redefine a new sequence: $n_h'=n_h$ if
$h\le i$ or $h>j$, while $n_h'=xn_h$ for $i<h\le j$. If $x=1$, then we
obtain the same sequence. Note that
$\frac{n_h'}{n_{h+1}'}=\frac{n_h}{n_{h+1}}$ for any $h$ except $h=i$
or $h=j$. So, we can slightly perturb $x$ either way, without
violating~\req{n}. The right-hand side of~\req{f}, as a function of
$x$ in a small neighborhood of $x=1$, is of the form $ax+\frac bx+c$
with $a,b>0$. But this function is strictly convex, so it cannot
attain its maximum at $x=1$, a contradiction.}
 
Let $t$ be the number of times we have $n_{i-1}=2n_i$. The
bound~\req{f} reads
 \beq{st}
 f(n,l)-  2 \log_2 n \le 2t + (s-t)\, \frac{l-1}2.
 \eeq
 Given that $2^t(\frac{l-1}2)^{s-t-1}\le n$, the right hand side
of~\req{st} is maximized for $t=O(\log l)$ and $s=(1+o(1))\, \frac{\ln
n}{\ln l}$, implying the upper bound~\req{MaxDeg} by~\req{DTf}.

Let us prove the lower bound. Let $k=\floor{l/2}$. Define
$G_0=K_{1,l-1}$ and $G_0'=K_{1,l-2}$. Let $r_0\in V(G_0)$, $r_0'\in
V(G_0')$ be their roots. Define inductively on $i$ the following
graphs. $G_{i}$ is obtained by taking $k$ copies of $G_{i-1}$ and
$k-1$ copies of $G_{i-1}'$, pairwise vertex-disjoint, plus the root
$r_i$ connected to the root of each copy of $G_{i-1}$ and
$G_{i-1}'$. We have $d(r_i)\le l-1$.  The graph $G_{i}'$ is defined in
a similar way except that we take $k-1$ copies of $G_{i-1}$ and $k$ copies
of $G_{i-1}'$. Let $i$ be the largest index such that
$\max(v(G_i),v(G_i'))\le n$.

Let us disregard all roots, i.e.,\ view $G_j$ and $G_j'$ as usual
(uncolored) graphs. Note that the trees $G_i$ and $G_i'$ are
non-isomorphic as for every $j$ we can identify the level-$j$ roots as
the vertices at distance $j+1$ from some leaf.

Define $g_j=(k-1)j+l-2$, $j\in[0,i]$. Let us show by induction on $j$
that Duplicator can survive at least $g_j$ rounds in the
$(G_j,G_j')$-game. This is clearly true for $j=0$. Let $j\ge 1$. If
Spoiler claims one of $r_j,r_j'$ then Duplicator selects the other. If
Spoiler selects a vertex in a graph from the ``previous'' level, for
example $F\subset G_j$ with $F\cong G_{j-1}'$, then Duplicator
chooses an $F'\subset G_i'$, $F'\cong G_{j-1}'$ and keeps the
isomorphism between $F$ and $F'$. So any moves of Spoiler inside
$V(F)\cup V(F')$ will be useless and we can ignore $F$ and $F'$. Thus
it takes Spoiler at least $k-1$ moves before we are down to the pair
$(G_{j-1},G_{j-1}')$, which proves the claim.

Thus we have $D(G_i) \ge D(G_i,G_i') \ge g_i=(\frac12+o(1))\,
\frac{l\ln n}{\ln l}$, finishing the proof.\epf

\brm Verbitsky~\cite{verbitsky:04} proposed a different argument
to estimate $\DT(n,l)$ which gives a weaker  bounds
than those in Theorem~\ref{th:MaxDeg} but can be applied
to other classes of graphs with small separators.

Let us study $\DT(n,l)$ for other $l,n$. The methods have much in
common with the proof of Theorem~\ref{th:MaxDeg} so our explanations
are shorter.

\bth{l=n^C} Let an integer $t\ge1$ be fixed. Suppose that
$l,n\to\infty$ so that $n\ge l^t$ but $n=o(l^{t+1})$. Then
$\DT(n,l)=(\frac{t+1}2+o(1))\, l$. In fact, the lower bound can be
achieved by uncolored trees.\end{theorem}
 \bpf The lower bound is proved by the induction on $t$. If $t=1$,
take $T_1=K_{1,l-2}$. One needs at least $l-1$ moves to distinguish it
from $T_1'=K_{1,l-1}$. Let $a=\floor{l/2}$ and $b=\ceil{l/2}$. Suppose
we have already constructed $T_{t-1}$ and $T_{t-1}'$, rooted trees
with $\le l^{t-1}$ vertices such that the root has degree at most
$l-1$. To construct $T_t$ take $a$ copies of $T_{t-1}$ and $b-1$
copies of $T_{t-1}'$ and connect them to the common root. For $T_t'$
we take $a-1$ and $b$ copies respectively. The degree of the main root
is $a+b-1= l-1$ as required. The order of $T_t$ is at most
$(a+b-1)l^{t-1}+1\le l^t$. Also, Spoiler needs at least $a$ moves
before reducing the game to $(T_{t-1},T_{t-1}')$ (while, for $t=1$,
$l$ moves are needed to finish the game), giving the required bound.

Let us turn to the upper bound. Spoiler uses the same strategy as
before. Namely, he chooses a median $x\in T$ and of two possible
multiplicities, summing up to $l$, chooses the smaller. Let
$m_1+1,m_2+1,\dots,m_k+1$ be the number of moves per each selected
median. We have $n\ge \prod_{i=1}^k m_i$.  Also, we have $k\le
\log_2n$ because we always choose a median. Given these restrictions,
the inequalities $m_i\le l/2$, $i\in[k-1]$, and $m_k\le l-1$, the sum
$\sum_{i=1}^k m_i$ is maximized if $m_k=l-1$ and as many as possible
$m_j=l/2$ are maximum possible. We thus factor out $l/2$ at most $t-1$
times until the remaining terms have the product (and so the sum)
$o(l)$. Thus,
 $$
 \sum_{i=1}^k (m_i+1)\le \log_2n+\sum_{i=1}^km_i\le l+\frac{(t-1)l}2+o(l),
 $$
 completing the proof.\epf

Theorems~\ref{th:MaxDeg} and~\ref{th:l=n^C} do not cover all the
possibilities for $n,l$. The asymptotic computation in the remaining
cases seems rather messy. However, the order of magnitude of
$\DT(n,l)$ is easy to compute with what we already have. Namely,
Theorem~\ref{th:l=n^C} implies that for $l=\Theta(n^t)$ with fixed
$t\in \I N$ we have $\DT(n,l)=\Theta(l)$. Also, if $l\ge 2$ is
constant, then $\DT(n,l)=\Theta(\ln n)$, where the lower bound follows
from considering the order-$n$ path and the upper bound is obtained by
using the method of Theorem~\ref{th:MaxDeg}.

\section{The Giant Component}\label{giant}

Let $c>1$ be a constant, $p=\frac cn$, and $G$ be the giant component
of a random graph $\C G(n,p)$. 
 \comment{
 Kim, Pikhurko, Spencer and
Verbitsky~\cite{kim+pikhurko+spencer+verbitsky:03rsa} conjectured that
whp $D(G)=O(\ln n)$.
 }%
 Here we show the following result.

 \bth{giant} Let $c>1$ be a constant, $p=c/n$, and $G$ be the giant
component of $\C G(n,p)$. Then whp
 \beq{giant}
 D(G)=\Theta\left(\frac{\ln n}{\ln \ln n}\right)
 \end{equation}
 \end{theorem}

This result allows us to conclude that for any $p=O(n^{-1})$ a
random graph $H\in \C G(n,p)$ satisfies whp
 \beq{d/n}
 D(H)=(\me^{-np}+o(1))\, n.
 \end{equation}
 The proof is an easy modification of that
in~\cite{kim+pikhurko+spencer+verbitsky:03rsa} where the validity of
\req{d/n} was established for $p\le (1.19...+o(1))\, n^{-1}$. The
lower bound in~\req{d/n} comes from considering the graph $H'$
obtained from $H$ by adding an isolated vertex (and noting that whp
$H$ has $(\me^{-np}+o(1))\, n$ isolated vertices). The method
in~\cite{kim+pikhurko+spencer+verbitsky:03rsa} shows that the upper
bound~\req{d/n} can fail only if $D(G)>(\me^{-np}+o(1))\, n$, where
$G$ is the giant component of $H$. (And $p/n\approx 1.19...$ is the
moment when $v(G)\approx \me^{-np}$.)

\subsection{Upper Bound}

The structure of the giant component is often characterized using its
core and kernel (e.g., see Janson, \L uczak, and
Ruci\'nski~\cite[Section~5]{janson+luczak+rucinski:rg}).  We follow this
approach in the proof of the upper bound in \req{giant}.  Thus, we
first bound $D(G)$ from above for a graph $G$ with small diameter
whose kernel fulfills some ``sparsness'' conditions. Then, we show
that these conditions hold whp for the kernel of the giant component
of a random graph.

\subsubsection{Bounding $D(G)$ Using the Kernel of $G$}\label{DKernel}

The \emph{core} $C$ of a graph $G$ is obtained by removing,
consecutively and as long as possible, vertices of degree at most
$1$. If $G$ is not a forest, then $C$ is non-empty and $\delta(C)\ge 2$.

First we need an auxiliary lemma which is easily proved, similarly to
the auxiliary lemmas in Section~\ref{general}, by the path-halving
argument.

\blm{cycle} Let $G,G'$ be graphs. Suppose $x\in V(G)$ and $x'\in
V(G')$ have been selected such that $G$ contains some cycle $P\ni x$
of length at most $k$ while $G'$ does not. Then Spoiler can win in at
most $\log_2 k+O(1)$ moves, playing all time inside
$G$.\qed\end{lemma}

\blm{DCore} Let $G,G'$ be graphs and $C,C'$ be their cores. If
Duplicator does not preserve the core, then Spoiler can win in at most
$\log_2d+O(1)$ extra moves, where $d$ is the diameter of
$G$.\end{lemma}
 \bpf Assume that $\diam(G')=\diam(G)$ for otherwise we are easily
done. Suppose that, for example, some vertices $x\in C$ and $x'\in \O
{C'}$ have been selected.

If $x$ lies on a cycle $C_1\subset C$, then we can find such a cycle
of length at most $2d+1$. Of course, $G'$ cannot have a cycle
containing $x'$, so Spoiler wins by Lemma~\ref{lm:cycle} in
$\log_2(2d+1)+O(1)$ moves, as required.

Suppose that $x$ does not belong to a cycle. Then $G$ contains two
vertex-disjoint cycles $C_1,C_2$ connected by a path $P$ containing
$x$. Choose such a configuration which minimizes the length of $P\ni
x$. Then the length of $P$ is at most $d$. Spoiler selects the branching
vertices $y_1\in V(C_1)\cap V(P)$ and $y_2\in V(C_2)\cap V(P)$. If
some Duplicator's reply $y_i'$ is not on a cycle, we done again by
Lemma~\ref{lm:cycle}. So assume there are cycles $C_i'\ni y_i'$. In
$G$ we have
 \beq{dist}
 \dist(y_1,y_2)= \dist(y_1,x) + \dist(y_2,x).
 \eeq
 As $x'\not\in C'$, any shortest $x'y_1'$-path and $x'y_2'$-path enter
$x'$ via the same edge $\{x',z'\}$. But then
 \beq{distp}
 \dist(y_1',y_2')\le \dist(y_1',z')+\dist(y_2',z')= \dist(y_1',x') +
\dist(y_2',x')-2.
 \eeq
 By~\req{dist} and~\req{distp}, the distances between $x,y_1,y_2$
cannot be all equal to the distances between $x',y_1',y_2'$. Spoiler
can demonstrate this in at most $\log_2 (\dist(y_1,y_2)) +O(1)$, as
required.\epf

In order to state our upper bound on $D(G)$ we have to define a number
of parameters of $G$. In outline, we try to show that any distict
$x,y\in V(C)$ can be distinguished by Spoiler reasonably fast. This
would mean that each vertex of $C$ can be identified by a first order
formula of small depth.
Note that  $G$ can be decomposed into the core and
a number of trees $T_x$, $x\in V(C)$, rooted at vertices of $C$.
Thus, by specifying which pairs of vertices of
$C$ are connected and describing each $T_x$, $x\in V(C)$, we
completely define $G$.  However, we have one unpleasant difficulty
that not all pairs of points of $C$ can be distinguished from one
another. For example, we may have a pendant triangle on $\{x,y,z\}$
with $d(x)=d(y)=2$, in which case the vertices $x$ and $y$ are
indistinguishable. However, we will show that whp we can
distinguish any two vertices of degree $3$ or more in $C$, which
suffices for our purposes.

Let us give all the details. For $x\in V(C)$, let $T_x\subset G$
denote the tree rooted at $x$, i.e., $T_x$ is a component containing
$x$ in the forest obtained from $G$ by removing all edges of $C$. Let
 $$
 t=\max\{D(T_x)\mid x\in V(C)\},
 $$
 where $D(T_x)$ is taken with respect to the class of graphs with one
root.

Let the \emph{kernel} $K$ of $G$ be obtained from $C$ by the
\emph{serial reduction} where we repeat as long as possible the
following step: if $C$ contains a vertex $x$ of degree $2$, then
remove $x$ from $V(C)$ but add the edge $\{y,z\}$ to $E(C)$ where
$y,z$ are the two neighbors of $x$. Note that $K$ may contain loops
and multiple edges. We agree that each loop contributes $2$ to the
degree. Then we have $\delta(K)\ge 3$.

Let $u=\Delta(G)$ and $d$ be the diameter of $G$. It follows that
each edge of $K$ corresponds to the path $P$ in $C$ of length at most
$2d$.
 \comment{(For otherwise any two vertices of $P$ at distance $d+1$
contradict the definition of $d$.)}%

Let $l$ be an integer such that every set of $v\le 6 l$ vertices of
$K$ spans at most $v$ edges in $K$. (Roughly speaking, we do not have
two short cycles close together.)

For $\{x,y\}\in E(K)$ let $A_{x,y}$ be the set of vertices obtained by
doing breadth first search in $K-x$ starting with $y$ until the
process dies or, after we have added a whole level, we reach at least
$k=2^{l-2}$ vertices. Let $K_{x,y}=K[A_{x,y}\cup \{x\}]$.

The \emph{height} of $z\in V(K_{x,y})$ is the distance in $K-x$ between
$z$ and $y$. It is easy to deduce from the condition on short cycles
that each $K_{x,y}\subset K-x$ has at most one cycle and the maximum
height is at most $l$. In fact, the process dies only in the case if
$y$ is an isolated loop in $K-x$.  For $xy\in E(K)$ let $G_{x,y}$ be a
subgraph of $G$ corresponding to $K_{x,y}$. We view $K_{x,y}$ and
$G_{x,y}$ as having two special \emph{roots} $x$ and $y$.

Here is another assumption about $G$ and $l$ we make. Suppose that for
any $xx',yy'\in E(K)$ if $K_{x,x'}$ and $K_{y,y'}$ have both order at
least $k$ and $A_{x,x'}\cap A_{y,y'}=\emptyset$, then the rooted graphs
$G_{x,x}$ and $G_{y,y'}$ are not isomorphic. Let
 \begin{eqnarray}
 b_0&=&\frac{l(\ln u+\ln \ln n + l)}{\ln l} +2u+\log_2d,\label{eq:b0}\\
 b&=& b_0 + t+ u +2\log_2d.\label{eq:b}
 \end{eqnarray}

\blm{a} Under the above assumptions on $G$, we have $D(G)\le b+O(1)$. 
\end{lemma}
 \bpf
 Let $G'\not\cong G$. Let $C',K'$ be its core and kernel. We can
assume that $\Delta(G')=u$ and its diameter is $d$ for otherwise
Spoiler easily wins in $u+2$ or $\log_2d+O(1)$ moves.

By Lemma~\ref{lm:DCore} it is enough to show that Spoiler can win the
Ehrenfeucht $(G,G')$-game in at most $b-\log_2d+O(1)$ moves provided
Duplicator always respects $C$ and $K$. Call this game $\C C$.

Color $V(K)\cup E(K)$ and $V(C)$ by the isomorphism type of the
subgraphs of $G$ which sit on a vertex/edge. We have a slight problem
with the edges of $K$ as the color of an unordered edge may depend in
which direction we traverse it. So, more precisely, every edge of $K$
is considered as a pair of ordered edges each getting its own
color. Do the same in $G'$. As $G\not\cong G'$, the obtained colored
digraphs $K$ and $K'$ cannot be isomorphic. Call the corresponding
digraph game $\C K$.

\claim1 If Spoiler can win the game $\C K$ in $m$ moves, then he can
win $\C C$ in at most $m+t+u+\log_2d+O(1)$ moves.
 \bcpf We can assume that each edge of $K'$ corresponds to a path in
$G'$ of length at most $2d+1$: otherwise Spoiler selects a vertex of
$C'$ at the $C'$-distance at least $d+1$ from any vertex of $K'$ and
wins in $\log_2d+ O(1)$ moves.

Spoiler plays according to his $\C K$-strategy by making moves inside
$V(K)\subset V(G)$ or $V(K')\subset V(G')$. Duplicator's reply are
inside $V(K')$, so they correspond to replies in the $\C K$-game. In
at most $m$ moves, Spoiler can achieve that the set of colored edges
between some selected vertices $x,y\in K$ and $x',y'\in K'$ are
different. (Or loops if $x=y$.)

In at most $u+1$ moves, Spoiler can either win or select a vertex $z$
inside a colored $xy$-path $P$ (an edge of $K$) such that $z'$ either
is not inside an $x'y'$-path (an edge of $K'$) or its path $P'\ni z'$
has a different coloring from $P$. In the former case, Spoiler wins by
Lemma~\ref{lm:path}: in $G$ there is an $xy$-path containing $z$ and
no vertex from $K$.

Consider the latter case. Assume that $|P|=|P'|$, for otherwise we are
done by Lemma~\ref{lm:path}.  Spoiler selects $w\in P$ such that
for the vertex $w'\in P'$ with $\dist_P(w,x)=\dist_{P'}(w',x')$ we
have $T_w\not\cong T'_{w'}$. If Duplicator does not reply with $w'$,
then she has violated distances. Otherwise Spoiler needs at most $t$
extra moves to win the game $\C T$ on $(T_w,T'_{w'})$ (and at most
$\log_2d+O(1)$ extra moves to catch Duplicator if she does not
respect $\C T$).\ecpf

It remains to bound $D(K)$, the colored digraph version. This requires
a few preliminary results.

\claim2 For any $\{x,x'\}\in K$ we have $D(K_{x,x})\le b_0+O(1)$ in
the class of colored digraphs with two roots, where $b_0$ is defined
by~\req{b0}.\medskip
 \bcpf Let $T=K_{x,x}$ and $T'\not\cong T$. If $T$ is a tree, then we
just apply a version of Theorem~\ref{th:MaxDeg} using the order
($\le\! u 2^{l}$) and maximum degree ($\le\!  u$). Otherwise, Spoiler
first selects a vertex $z\in T$ which lies on the (unique) cycle. We
have at most $u-1$ components in $T-z$, viewing each as a colored tree
where one extra color marks the neighbors of $z$. As $T\not\cong T'$,
in at most $u+1$ moves we can restrict our game to one of the
components. (If Duplicator does not respect components, she loses in
at most $\log_2 d +O(1)$ moves.)  Now, one of the graphs is a colored
tree, and Theorem~\ref{th:MaxDeg} applies.\ecpf

\claim3 For every two distinct vertices $x,y\in V(K)$ there is a first
order formula $\Phi_{x,y}(z)$ with one free variable and quantifier
rank at most $b_0+\log_2d+O(1)$ such that $G\models \Phi_{x,y}(x)$ and
$G\not\models \Phi_{x,y}(y)$. (Note that we have to find $\Phi_{x,y}$
for $x,y$ in the kernel only, but we evaluate $\Phi_{x,y}$ with
respect to $G$.)\medskip
 \bcpf To prove the existence of $\Phi_{x,y}$ we have to describe
Spoiler's strategy, where he has to distinguish $(G,x)$ and $(G,y)$
for given distinct $x,y\in K$. 

If the multiset of isomorphism classes $K_{x,x'}$, over $\{x,x'\}\in
E(K)$ is not equal to the multiset $\{ K_{y,y'}\mid \{y,y'\}\in
E(K)\}$, then we are done by Claim~2. So let us assume that these
multisets are equal. 

Note that an isomorphism $K_{x,x'}\cong K_{y,y'}$ implies an
isomorphism $G_{x,x'}\cong G_{y,y'}$.  Also, by our assumption on $l$,
the isomorphism $G_{x,x'}\cong G_{y,y'}$ implies that $V(K_{x,x'})\cap
V(K_{y,y'})\not=\emptyset$.

At most one neighbor of $x$ can be an isolated loop for otherwise, we
get 3 vertices spanning 4 edges. The same holds for $y$.  As the
height of any $K_{a,b}$ is at most $l$, we conclude that
$\dist_K(x,y)\le 2l$. A moment's thought reveals that there must be a
cycle of length at most $4l$ containing both $x$ and $y$. But this cycle
rules out the possibility of a loop adjacent to $x$ or to $y$. Thus,
in order to exclude $2$ short cycles in $K$ close to each other, it
must be the case that $\dist(x,y)\le l-1$ and
$d_K(x)=d_K(y)=3$. Moreover, let $x_1,x_2,x_3$ and $y_1,y_2,y_3$ be
the neighbors of $x$ and $y$ such that $G_{x,x_i}\cong G_{y,y_i}$;
then (up to a relabeling of indices), we have the following paths
between $x$ and $y$: either $(x,x_1,\dots,y_1,y)$ and
$(x,x_2,\dots,y_3,y)$ or  $(x,x_1,\dots,y_3,y)$ and
$(x,x_2,\dots,y_1,y)$

Now, $K_{x,x_3}$ is not isomorphic to $K_{x,x_1}$ nor to $K_{x,x_2}$
by the vertex-disjointness. (Note that it is not excluded that
$K_{x,x_1}\cong K_{x,x_2}$: they may intersect, for example, in $y$.)

But then $z=x$ is different from $z=y$ in the following respect: the
(unique) short cycle of $K$ containing $z$ has its two edges entering
$z$ from subgraphs isomorphic to $K_{x,x_1}$ and $K_{x,x_2}$ (while
for $z=y$ the corresponding subgraphs are isomorphic to $K_{x,x_1}$
and $K_{x,x_3}$).

This can be used by Spoiler as follows. Spoiler selects $x_1,x_2$. If
Duplicator replies with $y_3$, then Spoiler can use Claims~2 and~3
because $K_{y,y_3}$ is not-isomorphic to $K_{x,x_1}$ nor to
$K_{x,x_2}$. Otherwise, the edge $\{x,x_2\}$ is on a short cycle while
$\{y,y_2\}$ is not. Spoiler uses Lemma~\ref{lm:cycle}.\ecpf

By Lemma~\ref{lm:DCore} we can find $\Phi_K(x)$, a formula of rank at
most $\log_2d+O(1)$ which, with respect to $G$, evaluates to $1$ for
all $x\in V(K)$ and to $0$ otherwise. More precisely,
Lemma~\ref{lm:DCore} gives a formula $\Phi_C(x)$ testing for $x\in
V(C)$. But $V(K)\subset V(C)$ are precisely the vertices of degree at
least $3$ in $C$.
 \comment{So we can take
 $$
 \Phi_K(x)= \Phi_C(x) \wedge \exists_{x_1,x_2,x_1} \left(
 \Phi_C(x_1)\wedge  \Phi_C(x_2)\wedge  \Phi_C(x_3)\wedge x\sim
 x_1\wedge x\sim x_2\wedge x\sim x_3\wedge_{i\not= j} x_i\not=x_j\right).
 $$
 }

Now, as it is easy to see, for any $x\in K$ the formula
 \beq{Phi}
 \Phi_x(v):= \Phi_K(v) \wedge \bigwedge_{y\in V(K)\setminus \{x\}}
\Phi_{x,y}(v) 
 \eeq
 identifies uniquely $x$ and has rank at most
$\log_2d+b_0+ O(1)$. 

Take $x\in V(K)$. If there is no $x'\in V(K')$ such that $G'\models
\Phi_{x}(x')$, then Spoiler selects $x$. Whatever Duplicator's reply
$x'$ is, it evaluates differently from $x$ on $\Phi_{x}$. Spoiler can
now win in at most $D(\Phi_{x})$ moves, as required. If there are two
distinct $y',z'\in K'$ such that $G'\models \Phi_{x}(y')$ and
$G'\models \Phi_{x}(z')$, then Spoiler selects both $y'$ and $z'$. At
least one of Duplicator's replies is not equal to $x$, say,
$y\not=x$. Again, the selected vertices $y\in V(K)$ and $y'\in V(K')$
are distinguished by $\Phi_x$, so Spoiler can win in at most extra
$D(\Phi_x)$ moves. 

Therefore, let us assume that for every $x\in V(K)$ there is the
unique vertex $x'=\phi(x)\in V(K')$ such that $G'\models
\Phi_x(x')$. Clearly, $\phi$ is injective. Furthermore, $\phi$ is
surjective for if $x'\not\in \phi(V(K))$, then Spoiler wins by
selecting $x'\in V(K')$ and then using $\Phi_x$, where $x\in V(K)$ is
Duplicator's reply. Moreover, we can assume that Duplicator always
respects~$\phi$ for otherwise Spoiler wins in at most
$\log_2d+b_0+O(1)$ extra moves.

As $K\not\cong K'$, Spoiler can select $x,y\in V(K)$ such that the
multisets of colored paths (or loops if $x=y$) between $x$ and $y$ and
between $x'=\phi(x)$ and $y'=\phi(y)$ are distinct. Again, this means
that some colored path has different multiplicities and Spoiler can
highlight this in at most $u+1$ moves. Then in at most $\log_2l+O(1)$
moves he can ensure that some vertices $z\in V(K)$ and $z'\in V(K')$
are selected such that the removed trees $T_z$ and $T_{z'}$ rooted at
$z$ and $z'$ are not isomorphic, compare with
Lemma~\ref{lm:path}.

Now, by the definition of $t$, at most $t$ moves are enough to
distinguish $T_z$ from $T_{z'}'$ (plus possible $\log_2 d +O(1)$ moves
to catch Duplicator if she replies outside $V(T_z)\cup V(T_{z'})$). 

This completes the proof of Lemma~\ref{lm:a}.\epf

\subsubsection{Probabilistic Part}

Here we estimate the parameters from the previous section. As before,
let $G$ be the giant component of $\C G(n,\frac cn)$, let $C$ be its core,
etc. 

It is well-known that whp $u=O(\frac{\ln n}{\ln\ln n})$ and
$d=O(\ln n)$. \comment{Reference???}

\blm{Shaved} Whp every edge of $K$ corresponds to at most $O(\ln n)$
vertices of $G$. Similarly, for any $x\in V(C)$ we have
$v(T_x)=O(\ln n)$.\end{lemma}
 \bpf
 The expected number of $K$-edges, each corresponding to precisely $i=O(\ln n)$
vertices in $G$ is at most
 $$
 \binom{n}{i}\binom{i}{2} p^{i-1} i^{i-2} (1-p)^{(i-2)(n-i)} \le 
n i^2\left(\frac{\me c}{\me^c}\right)^{i}.
 $$
 But $\me c< \me^c$ for $c>1$, so if $i$ is large enough, $i>M\ln n$,
then the expectation is $o(n^{-3})$.

Similarly, the expected number of vertices $x$ with $v(T_x)=i=O(\ln n)$ is at most
$$n\binom{n-1}{i-1}p^{i-1}i^{i-2}(1-p)^{(i-1)(n-i)}\leq 2n i\left(\frac{\me c}{\me^c}\right)^{i}.$$
\epf

In particular, our results from Section~\ref{general} imply that whp
$t=O(\frac{\ln n}{\ln \ln n})$.

Let, for example, $l=2\ln \ln n$. Thus $k/\ln n\to\infty$, where
$k=2^{l-2}$. It remains to prove that this choice of $l$ satisfies all
the assumptions.

\blm{ShortCycle} Whp any set of $s\le 6l$ vertices of $K$ spans at
most $s$ edges.\end{lemma}
 \bpf A moment's thought reveals that it is enough to consider
sets spanning connected subgraphs only.

Let $L=M\ln n$ be given by Lemma~\ref{lm:Shaved}.  The probability
that there is a set $S$ such that $|S|=s\leq 6l$ and $K[S]$ is a
connected graph with at least $s+1$ edges is at most
 \begin{align*}
&o(1)+\sum_{s=4}^{6l}\binom{n}{s}\, s^{s-2}\, {s\choose 2}^2\sum_{0\leq
\ell_1,\ldots,\ell_{s+1}\leq L}
\prod_{i=1}^{s+1}\binom{n}{\ell_i}(\ell_i+2)^{\ell_i}p^{\ell_i+1}(1-p)^{\ell_i(n-\ell_i-2)}\\
 &\leq o(1)+\sum_{s=4}^{6l}\bfrac{n\me}{s}^s s^{s+2}\sum_{0\leq \ell_1,\ldots,\ell_{s+1}\leq L}
\prod_{i=1}^{s+1}\left(\frac{c\me^2}{n}\left(\frac{\me
  c}{\me^c}\right)^{\ell_i}\right)\ \le\ o(1)+
\sum_{s=4}^{6l}\frac{(O(1))^s}{n}\ =\ o(1).
\end{align*}
 The lemma is proved.\epf

\blm{Kab} Whp $K$ does not contain four vertices $x,x',y,y'$ such that
$xx',yy'\in E(K)$, $v(K_{x,x'})\ge k$, $A_{x,x'}\cap
A_{y,y'}=\emptyset$, and $G_{x,x'}\cong G_{y,y'}$.\end{lemma}
 \bpf Given $c$, choose the following constants in this order: small
$\e_1>0$, large $M_1$, large $M_2$, small $\e_2>0$, and large $M_3$. 

Consider breadth-first search in $G-x$ starting with $x'$. Let
$L_1=\{x'\}$, $L_2$, $L_3$, etc., be the levels. Let $T_i=\{x\}\cup
(\cup_{j=1}^i L_i)$. Let $s$ be the smallest index such that $|T_s|\ge
M_2\ln n$.

Chernoff's bound implies that the probability of $|T_s|> 2cM_2 \ln n$
is $o(n^{-2})$. Indeed, this is at most the probability that the
binomial random variable with parameters $(n, \frac cn \times M_2\ln n)$
exceeds $2cM_2\ln n$.

Similarly, with probability $1-o(n^{-3})$ we have $|L_{i+1}|=(c\pm
\e_2)|L_i|$ provided $i\ge s$ and $|T_i|=o(n)$.  Hence, we
see that from the first time we reach $2M_2\ln n$ vertices, the levels
increase proportionally with the coefficient close to $c$ for further
$\Theta(\ln n)$ steps.

Take some $i$ with $|T_i|=O(\ln n)$. The sizes of the first $\Theta(\ln
n)$ levels of the breadth-first search from the vertices of $L_i$ can
be bounded from below by independent branching processes with the number of
children having the Poisson distribution with mean $c-\e_2$. Indeed,
for every active vertex $v$ choose a pool $P$ of
$\ceil{(1-\frac{\e_2}c)n}$ available vertices and let $v$ choose its
neighbors from $P$, each with probability $c/n$. (The edges between
$v$ and $\O P$ are ignored.) If $v$ claimed $r$ neighbors, then, when
we take the next active vertex $u$, we add extra $r$ vertices to the
pool, so that its size remains constant.

With positive probability $p_1$ the ideal branching process survives
infinitely long; in fact, $p_1$ is the positive root of
$1-p_1=\me^{-cp_1}$. Let
 $$
 p_2=\max_{j\ge 0} \frac{c^j\me^{-c}}{j!} <1.
 $$
 The numbers $p_1>0$ and $p_2<1$ are constants (depending on $c$
only).

Take the smallest $i$ such that $|T_i|\ge 2cM_3\ln n$. The
breadth-first search inside $G$ goes on for at least $M_1$ further
rounds (after the $i$-th round) before we reach a vertex outside
$G_{x,x'}$. We know that $|L_i|\ge (\frac{c-1}c-\e_1)\,|T_i|$ because
the levels grow proportionally from the $s$-th level. Let $Z$ consist
of the vertices of $L_i$ for which the search process in $G-x$ goes on
for at least $M_1$ further levels before dying out. By Chernoff's
bound, with probability $1-o(n^{-2})$ we have $|Z|\ge \frac{p_1}2
|L_i|$.

Let us fix any $K_{x,x'}$ having all the above properties and compute
the expected number of copies of $K_{x,x'}$ in $G$. More precisely, we
compute the expected number of subgraphs of $G$ isomorphic to
$G[T_{i}]$ such that a specified $|Z|$-subset of the last level has
specified trees, each of height at least $M_1$, sitting on it. The
expected number of $G[T_i]$-subgraphs is at most
$n^{|T_i|}\,p_1^{|T_i|-1}$. This has to be multiplied by
 $$
 (p_2+o(1))^{M_1|Z|} \le p_2^{M_1(c-1)p_1\,|T_i|/4c}:
 $$
 because if we want to get a given height-$M_1$ tree, then at least
$t$ times we have to match the sum of degrees of a level, each
coincidence having probability at most $p_2+o(1)$. As the constant
$M_1$ can be arbitrarily large, we can make the total expectation
$o(n^{-2})$.

Markov's inequality implies the lemma.\epf

Finally, putting all together we deduce the upper bound of
Theorem~\ref{th:giant}.

\subsection{Lower Bound}

Let $l=(1-\e) \frac{\ln n}{\ln \ln n}$ for some $\e>0$. We claim that
whp the core $C$ has a vertex $i$ adjacent to at least $l$ leaves of
$G$. (Then we have $D(C)\ge l+1$: consider the graph obtained from $C$
by adding an extra leaf to $i$.)

Let us first prove this claim for the whole random graph $H\in \C
G(n,c/n)$ (rather than for the giant component $G\subset H$). For
$i\in [n]$ let $X_i$ be the event that the vertex $i$ is incident to
at least $l$ leaves. It is easy to estimate the expectation of
$X=\sum_{i=1}^n X_i$:
 \begin{eqnarray*}
 E(X) &=& n \binom{n-1}{ l} p^l (1-p)^{\binom{l}{ 2} + l(n-l)} +O(1)\times
 n\binom{n}{ l+1} p^{l+1}(1-p)^{(l+1)n}\\
 &=& (1+o(1)) \frac{nc^l\me^{-cl}}{l!}\ \to\ \infty.
 \end{eqnarray*}
 Also, for $i\not=j$,
 \begin{eqnarray*}
 E(X_i\wedge X_j) &=&(1+o(1))\,\binom{n-2}{ l} \binom{n-l-2}{ l}p^{2l}
 (1-p)^{\binom{2l}{ 2} +2l(n-2l-1)}\\
  &=& (1+o(1))\, E(X_i)E(X_j).
 \end{eqnarray*}
  The second moment method gives that $X$ is concentrated around its
mean. 

Now, let us reveal the vertex set $A$ of the $2$-core of the whole
graph $H$. When we expose the stars demonstrating $X_i=1$ one by one,
then for each $i$ the probability of $i\in A$ is $\frac{|A|}n+o(1)$.
The sharper results of {\L}uczak~\cite{luczak:91}
 \comment{Or Pittel~\cite{pittel:90}?}%
 imply that whp the core $C$ of the giant component has size
$\Theta(n)$. Hence, whp at least one vertex $i$ with $X_i=1$ belongs
to the $V(C)$, giving the required.

\section{Random Trees}\label{random}

We consider the probabilistic model $\C T(n)$, where a tree $T$ on the
vertex set $[n]$ is selected uniformly at random among all $n^{n-2}$
trees. In this section we prove that whp $D(T)$ is close to the
maximum degree of $T$.

\bth{RandomTree} Let $T\in\C T(n)$. Whp 
$D(T)=(1+o(1))\Delta(T)=(1+o(1))\frac{\ln
n}{\ln\ln n}$.
\end{theorem}

\newcommand{\Var}{{\textrm{Var}}}
\newcommand{\Chq}{{\textrm{Ch}}}
\newcommand{\rmv}{{\textrm{del}}}

Let $\cF(n,k)$ be a forest chosen uniformly at random from the family
of $\cF_{n,k}$ of all forests with the vertex set
$[n]$, which consist of $k$ trees rooted at vertices
$1,2,\dots,k$.  Note that a random tree $T\in \cT(n)$ can be
identified with $\cF(n,1)$.  We recall that $|\cF_{n,k}|=kn^{n-k-1}$,
see e.g.\ Stanley~\cite[Theorem~5.3.2]{stanley:ec}.  We start with the
following simple facts on $\cF(n,k)$.

\blm{forest} Let $k=k(n)\le \ln^4 n$.
\begin{enumerate}
\renewcommand{\labelenumi}{(\roman{enumi})} 
\item The expected number of vertices in all trees of $\cF(n,k)$,
except for the largest one, is $O(k\sqrt n)$.
\item The probability that $\cF(n,k)$ contains precisely $\ell$, $\ell=0,\dots,k-1$,
isolated vertices is given by $(1+O({k^2}/{n}))
\binom{k-1}\ell \me^{-\ell}(1-\me^{-1})^{k-\ell-1}$.
\item The probability that the roots of $\cF(n,k)$ have more than $k(1+1/\ln n)+2\ln^2 n$ 
neighbors combined is $o(n^{-3})$.
 \item The probability that $\ell$ given roots of $\cF(n,k)$ have
degree at least $s\ge 4$ each is bounded from above by $(2/(s-1)!)^\ell$
 \end{enumerate}
\end{lemma}

\bpf If $i\le n/2+1$, then the probability that a tree rooted at
a vertex $j=1,2,\dots,k$ in the forest $\cF(n,k)$ has precisely $i$
vertices is given by
 $$\binom {n-k}{i-1} i^{i-2} 
\frac{(k-1)(n-i)^{n-i-k}}{k n^{n-k-1}}=O(i^{-3/2})\,.$$
 Consequently, the expectation of the sum of the orders of all
components of $\cF(n,k)$ with at most $n/2+1$ vertices is $O(k \sqrt n)$.

In order to see (ii) note that from the generalized
inclusion-exclusion principle the stated probability equals
 \begin{equation}\label{eqf1}
\begin{aligned}
\sum_{i=\ell}^k&\binom i\ell(-1)^{i-\ell}\binom ki\frac{(k-i)(n-i)^{n-k-1}}{kn^{n-k-1}}\\
=&\Big(1+O\Big(\frac{k^2}{n}\Big)\Big)
\sum_{i=\ell}^k\frac{(k-1)!}{\ell!(i-\ell)!(k-1-i)!}(-1)^{i-\ell}\me^{-i}\\
=&\Big(1+O\Big(\frac{k^2}{n}\Big)\Big)
\binom{k-1}\ell \me^{-\ell}(1-\me^{-1})^{k-\ell-1}\,.
\end{aligned}
\end{equation}

For the probability that precisely $m$ ($\ge\! k$) vertices 
of $\cF(n,k)$ are adjacent to the roots, Stirling's formula gives
\begin{equation}\label{f1}
\binom{n-k}{m}k^m\frac{m\,(n-k)^{n-k-m-1}}{k\,n^{n-k-1}}
\le \Big(1+O\Big(\frac{k^2}n\Big)\Big)\Big(\frac{\me^ {1-k/m}k}{m}\Big)^{m}.
\end{equation}
For every $x$, $0<x<1$, we have $x\me^{1-x}\le \me^{-(1-x)^2/2}$,
so the above formula is bounded from above by $\exp(-\frac{(m-k)^2}{2m})$. 
Since   
$$\sum_{m\ge k(1+1/\ln n)+2\ln ^2n}\exp\Big(-\frac{(m-k)^2}{2m}\Big)=o(n^{-3})\,,$$
the assertion follows.

For $k=1$ the probability that a given root has
degree at least $s$ is bounded from above by
 $$\sum_{t\ge s}\binom{n-1}{t}\frac{t(n-1)^{n-t-2}}{n^{n-2}}\le
\sum_{t\ge s}\frac{1}{(t-1)!}\le \frac{2}{(s-1)!}\;.$$
 If we fix some $\ell\ge 2$ roots, then if we condition on the vertex
sets of the $\ell$ corresponding components, the obtained trees are
independent and uniformly distributed, implying the required bound by
the above calculation.
\epf

Using the above result one can estimate 
the number of vertices of $T\in \cT(n)$ with a
prescribed number of pendant neighbors.

\blm{vert} Let $X_{\ell,m}$ denote the number of vertices in $T\in
\cT(n)$ with precisely $\ell$ neighbors of degree one and $m$
neighbors of degree larger than one.  Let
 $$
 A\subseteq\{(\ell,m)\colon\; 0\le \ell\le \ln n, \quad 1\le m\le \ln n
\}\,,$$
 be a set of pairs of natural numbers and $X_A=\sum_{(\ell,m)\in A}
X_{\ell,m}$. Then, the expectation
 \begin{equation}\label{eqf2}
 E(X_A)=(1+o(1))\,n\sum_{(\ell,m)\in A}
\frac{\me^{-\ell-1}}{\ell!}\frac{(1-\me^{-1})^{m-1}}{(m-1)!}
\end{equation}
 and $E(X_A(X_A-1))=(1+o(1))\,(E(X_A))^2$.  \end{lemma}

\bpf 
Using Lemma~\ref{lm:forest}(ii) we get
\begin{equation*}
 E(X_A)=(1+o(1))n\sum_{(\ell,m)\in A}\binom{n-1}{m+\ell}\binom{m+\ell-1}\ell
\me^{-\ell} (1-\me^{-1})^{m-1}\frac{(m+\ell)(n-1)^{n-m-\ell-2}}{n^{n-2}}
\end{equation*}
 which gives (\ref{eqf2}). In order to count the expected number of pairs
of vertices with prescribed  neighborhoods one needs first to choose
 $\ell+m$ neighbors of a vertex  and then compute the expectation of
the number of vertices of a given neighborhood in the random forest 
$\cF(n,\ell+m)$ obtained in this way. However, the largest
tree of $\cF(n,\ell+m)$ has the expectation $n-O(\sqrt n \ln n)$ (Lemma~\ref{lm:forest}); 
one can easily observe that this fact implies  
that the expected number of vertices with a prescribed neighborhood in $\cF(n,\ell+m)$ 
 is $(1+o(1))\,E(X_A)$, and so $E(X_A(X_A-1))=(1+o(1))\,(E(X_A))^2$.  
\epf

As an easy corollary of the above result we get a lower bound for 
$D(\cT(n))$.

\bth{lower} Let $T\in\C T(n)$.
Whp $D(T)\ge (1-o(1))\Delta(T)=(1-o(1))\, \frac{\ln n}{\ln \ln n}$.
\end{theorem}

\bpf Since whp the maximum degree 
is $(1-o(1)){\ln n}/{\ln\ln n}$, in order to prove the assertion 
it is enough to show that whp $T$
contains a vertex $v$ with
 \begin{equation}\label{eqf3}
\ell_0=(1-o(1))\, \frac{\ln n}{\ln \ln n}
\end{equation}
neighbors of degree one; indeed, to characterize such a structure
Spoiler needs at least $\ell_0+1$ moves. Using Lemma~\ref{lm:vert}, we
infer that the for the number of vertices $X_{\ell}$ of $T$ with
exactly $\ell$ neighbors of degree $1$ we have
$E(X_\ell)=O(\me^{-\ell}n/\ell!)$. Thus, one can choose $\ell_0$ so that
(\ref{eqf3}) holds and $E(X_{\ell_0})\to\infty$.  Then, due to
Lemma~\ref{lm:vert}, $\Var(X_{\ell_0})=o((E(X_{\ell_0}))^2)$, 
and Chebyshev's inequality implies that whp $X_{\ell_0}>0$.\epf

Let us state  another simple consequence of Lemma~\ref{lm:forest}
which will be used in our proof of Theorem~\ref{th:RandomTree}. Here and below  $N_r(v)$
denotes the $r$-neighborhood of $v$, i.e., the set of all vertices of a graph which 
are at the distance $r$ from $v$, and $N_{\le r}(v)=\bigcup_{i=0}^r N_i(r)$.

\blm{largedegrees} Let $r_0=r_0(n)= \lceil 7 \ln n\rceil $. Then,
whp the following holds for every vertex $v$ of $T\in \cT(n)$:
\begin{enumerate}
\renewcommand{\labelenumi}{(\roman{enumi})} 
\item $|N_{\le r_0}(v)|\le 10^8 \ln^4n\;,$
\item  $N_{\le r_0}(v)$ contains fewer than  $\ln n/(\ln\ln n)^2$ vertices 
of degree larger than $(\ln\ln n)^5$.
\end{enumerate}
\end{lemma}

\bpf For $s\le r_0$ let $W_s=\cup_{i=0}^s N_i(v)$.
Note that,  conditioned on the structure of the
subtree of $T$ induced by $W_s$
for some $s\le r_0$,  the forest $T- W_{s-1}$
can be identified with the random forest on $n-|W_{s-1}|$
vertices, rooted at the set $W_s$. Thus,  
it follows from Lemma~\ref{lm:forest}(iii) that 
once for some $i$ we have $|N_i(v)|\ge 4 \ln ^3 n$ 
then $|N_{i+1}(v)|\le |N_i(v)|(1+2/\ln n)$,
 so that
$$|N_{\le r_0}(v)|\le 4 r_0\ln ^3n (1+2/\ln n)^{r_0}\le 10^8 \ln^4n\;.$$

In order to show (ii) note that (i) and  Lemma~\ref{lm:forest}(iv) 
imply  that the probability that, for some
vertex $v$, at least  $\ell=\lfloor \ln n/(\ln\ln n)^2\rfloor$ vertices of $N_{\le r_0}(v)$
 have degree larger than $m=(\ln\ln n)^5$ is bounded from above by 
$$n\binom {\ln^5 n}{\ell}\left(\frac{2}{(m-1)!}\right)^\ell
\le n\left(\frac{2\me \ln^5n}{\ell(m-1)!}\right)^\ell\le n\me^{-m\ell}=o(1).$$
 \comment{
 Here is a small hole: we know that the probability of having at least
$>m$ neighbors is at most $2/m!$ but why is the probability that $l$
given vertices each have degree $>m$ is at most $(2/(m-1)!)^l$?

Proof: we expose levels one by one. Once we have exposed a level, we
allow an adversary to choose any number of active vertices, provided
he does not choose more than $\ell$ vertices in total. Then adversary
succeeds (all his points have high degree) with probability at most
$(2/(m-1)!)^\ell$. 
 }%
 \epf

In our further argument we need some more definitions.  Let $T$ be a
tree and let $v$ be a vertex of $T$. For a vertex $w\in N_r(v)$ let
$P_{vw}$ denote the unique path connecting $v$ to $w$ (of length
$r$). Let the \emph{check} $\Chq(v;P_{vw})$ be the binary sequence
$b_0\cdots b_r$, in which, for $i=0,\dots, r$, $b_{i}$ is zero (resp.\
1) if the $i$-th vertex of $P_{vw}$ is adjacent (resp.\ not adjacent)
to a vertex of degree one.  Finally, the \emph{$r$-checkbook}
$\Chq_r(v)$ is the set
 $$
 \Chq_r(v)=\{\Chq(v;P_{vw})\colon w\in N_r\textrm{\ and }P_{vw} 
\textrm{ is a path of length $r$}\}.
$$
 Note that a checkbook  is not a multiset, i.e., a check from  
$\Chq_r(v)$ may correspond to more than one paths $P_{vw}$.

Our proof of the upper bound for $D(\cT(n))$ is based on the following 
fact. 

\bth{checks} Let $r_0=\lceil 7 \ln n\rceil$. 
Whp for each pair  $P_{vw}$, $P_{v'w'}$ of 
paths of length $r_0$ in $T\in \cT(n)$ which share at most one vertex, 
the checks  $\Chq(v;P_{vw})$ and $\Chq(v;P_{v'w'})$ 
are different.
 \end{theorem}

\bpf Let $C=\rmv(T)$ denote the tree obtained from $T$ by removing
all vertices of degree one.  From Lemma~\ref{lm:vert} it follows that
whp the tree $C$ has $(1-\me^{-1}-o(1))n$ vertices of which
 $$
 (1+o(1))\,n \sum_{\ell>0} \frac{\me^{-\ell-1}}{\ell!} =
(\exp(\me^{-1}-1)-\me^{-1}+o(1))\,n$$
 vertices have degree one and
 $$
 \alpha n = (1-\exp(\me^{-1}-1) +o(1))\, n.
 $$
 vertices have degree greater than one.

Moreover, among the set  $B$ of $(\me^{-1}+o(1))n$ vertices removed from $T$,
 $$
 (1+o(1))n\sum_{l=0}^\infty
\ell\frac{\me^{-\ell-1}}{\ell!}=(1+o(1))\exp(\me^{-1}-2)n\,$$
 were adjacent to vertices which became pendant in $C$. 
Let $B'$ denote  the set of the remaining 
 $$
 (\me^{-1}-\exp(\me^{-1}-2)+o(1))n=(\rho_0+o(1))n
 $$ 
vertices which are adjacent to vertices of degree at least two in $C$.
Note  that, given $C=\rmv(T)$, each attachment of 
vertices from $B\setminus B'$ to pendant vertices 
of $C$ such that each pendant vertex of $C$ get at least one vertex 
from $B\setminus B'$, as well as each attachment of vertices from $B'$ 
to vertices of degree at least two from $C$ is equally likely. 

Let $P_{vw}$, $P_{v'w'}$, be two paths of length $r_0$ in $T$
which share at most one vertex.  Clearly, each vertex of $P_{vw}$,
except, maybe, at most two vertices at each of the ends, belong to $C$ and
have in it at least two neighbors; the same is true for $P_{v'w'}$.
Since $(\rho_0+o(1))n$ vertices from $B'$ are attached to the $\alpha
n$ vertices of degree at least two in $C$ at random, the
probability that one such vertex gets no attachment is 
 $$
 p_0=(1+o(1))\, \left(1-\frac1{\alpha n}\right)^{\rho_0 n}= (1+o(1))\,
\me^{-\rho_0/\alpha} = 0.692...+o(1).
 $$
 Therefore, the probability that the checks $\Chq(v,P_{vw})$ and
$\Chq(w,P_{v'w'})$ are identical is bounded from above by
 $$
 \left(p_0^2+(1-p_0)^2 +o(1)\right)^{r_0}\le \me^{-3\ln n}=o(n^{-2})\,.
 $$
 
Since by Lemma~\ref{lm:largedegrees}(i) whp $T$ contains at most
$O(n\ln^4 n)$ checks of length $r_0$, the assertion follows.
 \epf

Now, let $r_0=\lceil 7 \ln n\rceil$. 
We call a tree $T$ on $n$ vertices \emph{typical} if:
\begin{itemize}
\item for each pair of paths $P_{vw}$, $P_{v'w'}$ of length 
$r_0$ which share at most one vertex, 
the checks  $\Chq(v;P_{vw})$, $\Chq(v;P_{v'w'})$ are different,
\item for the maximum degree $\Delta$  of $T$ we have
 $$\frac{\ln n}{2\ln\ln n}\le \Delta\le \frac{2\ln n}{\ln\ln n} \,,$$
\item $|N_{\le r_0}|\le 10^8\ln ^4 n$,  for every vertex $v$, 
\item for every vertex $v$  at most $\ln n/(\ln\ln n)^2$ vertices of
degree larger than $(\ln\ln n)^5$ lie within distance 
$r_0$ from $v$.  
\end{itemize}

\bth{upper} For a typical tree $T\in\C T(n)$ we have
$D(T)\le (1+o(1))\, \Delta$.  \end{theorem}

\bpf Let $T$ be a typical tree and $T'$ be any other graph which is
not isomorphic to $T$. We shall show that then Spoiler can win the
Ehrenfeucht game on $T$ and $T'$ in $(1+o(1))\Delta$ moves.

Let us call a vertex $v$ of a graph a \emph{yuppie}, if there are two
paths $P_{vw}$, $P_{vw'}$ of length $r_0$ starting at $v$ so that
$V(P_{vw})\cap V(P_{vw'})=\{v\}$.  Note that the set of all yuppies
$Y$ spans a subtree in $T$, call it $K$.

Our approach is similar to that for the giant component from
Section~\ref{giant}. 

Let us view $K$ as a colored graph where the color of a vertex $x$ is
the isomorphism type of the component of $T-(Y\setminus\{x\})$ rooted
at $x$. Let $Y'$ be the set of yuppies of $T'$, and let
$K'=T'[Y']$. We can assume that Duplicator preserves the subgraphs $K$
and $K'$, for otherwise Spoiler wins in extra $O(\ln \ln n)$ moves.

\claim1 Any distinct $v,v'\in K$ can be distinguished (with respect to
$G$) in $O(\ln\ln n)$ moves.\medskip
 
\bcpf Assume that the $r_0$-checkbooks of $v,v'$ are the same for
otherwise Spoiler wins in $\log_2(r_0)+O(1)$ moves. (Please note that
the checkbooks are viewed as sets, not as multisets, so the number of
moves does not depend on the degrees of $v$ and $v'$.)

Take a path $P_{vx}$ of length $r_0$, which shares with $P_{vv'}$
only vertex $v$. Spoiler selects $x$. Let Duplicator reply with
$x'$. Assume that $\Chq(w,P_{vx})=\Chq(v,P_{v'x'})$.  The path
$P_{v'x'}$ must intersect $P_{vx}$; thus $v\in P_{v'x'}$. Next,
Spoiler selects the $P_{vx}$-neighbor $y$ of $v$; Duplicator's reply
must be $y'\in P_{v'x'}$. 

Let $z\in T$ maximize $\dist(v,z)$ on the condition that
$\Chq(z)=\Chq(v)$ and $v$ lies between $y$ and $z$ in $T$. Define the
analogous vertex $z'$, replacing $v,y$ in the definition by
$v',y'$. We have $\dist(v,z)>\dist(v',z')$. Let Spoiler select
$w=z$.  If Duplicator's reply $w'$ satisfies $\Chq(w')\not\cong
\Chq(w)$, then Spoiler quickly wins. Otherwise, $\dist(v,w)>
\dist(v',w')$. Moreover, $\dist(v,w)\le 2r_0$ (because their
$r_0$-checkbooks are non-empty and equal). Spoiler wins in $\log_2
r_0+O(1)$ extra moves. The claim has been proved.\ecpf

Similarly to the argument surrounding~\req{Phi}, one can agrue that
for every vertex $x\in K$ there is a formula $\Phi_x(v)$ of rank
$O(\ln \ln n)$ identifying $x$ (with respect to $T$). Moreover, we can
assume that this gives us an isomorphism $\phi:K\to K'$ which is
respected by Duplicator.

As $T\not\cong T'$, there are two cases to consider.

\case1 There is $x\in K$ such that $T_x\not\cong T'_{x'}$, where
$x'=\phi(x)$ and $T_{x'}'$ is the component of $T'-(Y' \setminus\{x'\})$
rooted at $x'$.\medskip

Since each vertex of $T$ is within distance at most $r_0$ from some
yuppie, the tree $T_x$ has height at most $r_0$. If $T'_{x'}$ has a
path of length greater than $2r_0$ or a cycle, then Spoiler easily
wins, so assume that $T'$ is a tree. Now Spoiler should select all
vertices of $T_x$ which are of degree larger than $(\ln\ln n)^5$, say
$w_1,\dots,w_t$. Since $T$ is typical there are at most $\ln n/(\ln\ln
n)^2$ such vertices in $T_v$.  Suppose that, in responce to that,
Duplicator chooses vertices $w'_1,\dots,w'_s$ in $T'_{x'}$.  Then,
$T_v\setminus \{w_1,\dots,w_s\}$ splits into a number of trees $F_1,
\dots, F_u$, colored accordingly to their adjacencies to the
$w_i$'s. Now, for some $i$ the multisets of colored trees adjacent to
$w_i$ and $w_i'$ are different. Spoiler can highlight this by using at
most $\Delta(T)+1$ moves. Now Spoiler plays inside some $F_i$ the
strategy of Theorem~\ref{th:MaxDeg}. Note that $F_i$ has diameter
at most $2r_0$ and maximum degree at most $(\ln\ln n)^5$.\medskip

\case2 $T'$ is not connected.\medskip

As $K'\cong K$ is connected, there is a component $C'$ of $T'$ without
a yuppie. Spoiler chooses an $x'\in C'$. Now, any Duplicator's reply
$x$ is within distance $r_0$ from a yuppie, which is not true for
$x'$. Spoiler can win in $O(\ln \ln n)$ moves.\medskip

Consequently, for a typical tree $T$,
 $$
 D(T)\le \Delta(T)+\frac{\ln n}{(\ln\ln n)^2}+O((\ln\ln n)^6)\,, 
 $$ 
and the assertion follows.
\epf

\noindent {\it Proof of Theorem~\ref{th:RandomTree}.}
Theorem~\ref{th:RandomTree} is an immediate consequence
of Theorems~\ref{th:lower} and~\ref{th:upper} and the fact that, 
due to Lemmas~\ref{lm:forest} and~\ref{lm:largedegrees},
whp a random tree $T\in \cT(n)$ is typical.
\qed

\section{Restricting Alternations}

If Spoiler can win the Ehrenfeucht game, alternating between the
graphs $G$ and $G'$ at most $r$ times, then the corresponding sentence
has the \emph{alternation number} at most $r$, that is, any chain of
nested quantifiers has at most $r$ changes between $\exists$ and
$\forall$. (To make this well-defined, we assume that no quantifier is
within the range of a negation sign.) Let $D_r(G)$ be the smallest
depth of a sentence which defines $G$ and has the alternation number
at most $r$. It is not hard to see that $D_r(G)=\max\{D_r(G,G')\mid
G'\not\cong G\}$, where $D_r(G,G')$ may be defined as the smallest $k$
such that Spoiler can win $\ehr_k(G,G')$ with at most $r$
alternations. For small $r$, this is a considerable restriction on the
structure of the corresponding formulas, so let us investigate the
alternation number given by our strategies.

Let $\DT_r(n,l)$ be the maximum of $D_r(T)$ over all colored trees of
order at most $n$ and maximum degree at most $l$.

Unfortunately, in Theorem~\ref{th:MaxDeg} we have hardly any control
on the number of alternations. However, we can show that alternation
number $0$ suffices if we are happy to increase the upper bound by a
factor of $2$.

 \begin{lemma}\label{lem:treetree}
Let $T$ and $T'$ be colored trees. Suppose that $T\not\cong T'$,
where $\cong$ stands for the isomorphism relation for colored trees, i.e.,
the underlying (uncolored) trees of $T$ and $T'$ may be isomorphic.
Furthermore, assume that $v(T)\ge v(T')$ and denote $n=v(T)$.
Assume also that $\Delta(T)\le l$ and let
both $l$ and $\ln n/\ln l$ tend to the infinity. 
Then Spoiler can win the Ehrenfeucht game on $(T,T')$ in at most
 \beq{D1}
 (1+o(1)) \frac{l \ln n}{\ln l}.
 \eeq
moves playing all time in~$T$.
\end{lemma}

 \bpf
In the first move Spoiler selects a median $x\in T$; let
$x'$ be Duplicator's reply. 
If $d(x)>d(x')$, then
Spoiler wins in extra $l$ moves, which is negligible when compared
to~(\ref{eq:D1}). So, suppose that $d(x')\ge d(x)$.

Let $t=d(x)$ and $C_1,\dots,C_t$ be the (rooted) components of $T-x$
indexed so that $v(C_1)\ge v(C_2)\ge\ldots\ge v(C_t)$.  Referring to
the root of a component we mean the vertex of it which is adjacent to
$x$.  Spoiler starts selecting, one by one, the roots of
$C_1,C_2,\ldots$.  Duplicator is enforced to respond with roots of
distinct components of $T'-x'$. Spoiler keeps doing so until the
following situation occurs: he selects the root $y$ of a component
$C=C'_i$ while Duplicator selects the root $y'$ of a component $C'$
such that $v(C)\ge v(C')$ and $C\not\cong C'$ (as rooted trees). Such
a situation really must occur for some $i\le t$ due to the conditions
that $v(T)\ge v(T')$, $d(x)\le d(x')$, and $T\not\cong T'$.

We claim that if Spoiler selects a vertex $z$ inside $C$, then
Duplicator must reply with some $z'\in C'$ for otherwise Spoiler wins
in at most $\log_2 n$ moves. Indeed, suppose $z'\not\in C'$. Spoiler
selects $z_1$ which is a middle point of the $yz$-path. Whatever the
reply $z_1'$ is, the $z'z_1'$-path or $z_1'y'$-path contains the vertex
$x'$. Suppose it is the $z'z_1$-path. Then Spoiler halves the
$zz_1$-path. In at most $\log_2n$ times he wins.

Thus making $i+1\le t+1\le l+1$ steps, we have reduced the game to two
non-isomorphic (rooted) trees, $C$ and $C'$, with $v(C)\le
\min(\frac1i,\frac12)\, v(T)$.  In the game on $(C,C')$ Spoiler
applies the same strategy recursively.  Two ending conditions are
possible: the root of $C$ has strictly larger degree than the root of
$C'$ and Duplicator violates a color, the adjacency, or the equality
relation.  It is easy to argue, cf.\ the proof of
Theorem~\ref{th:MaxDeg}, that the worst case for us is when we have
$i=(1+o(1))\, l$ all the time, which gives the required
bound~(\ref{eq:D1}).
 \epf

\bth{DT0} Let both $l$ and $\ln n/\ln l$ tend to the
infinity. Then
 \beq{}
 \DT_0(n,l)\le (1+o(1)) \frac{l \ln n}{\ln l}.
 \eeq
 \end{theorem}

\bpf
 Let $T$ be a tree of order $n$ and maximum degree at most $l$
and let $G\not\cong T$. 
If $\Delta(T)\ne\Delta(G)$ then Spoiler wins the Ehrenfeucht game on
$(T,G)$ in at most $l+2$ moves
playing in the graph of the larger degree. We will therefore assume that
$T$ and $G$ have the same maximum degree not exceeding~$l$.

\case1 $G$ contains a cycle of length no more than $n+1$.\smallskip

Spoiler plays in $G$ proceeding as in the last paragraph of the proof
of Lemma~\ref{lm:Tree}.

\case2 $G$ is connected and has no cycle of length up to $n+1$.

If $v(G)\le n$, then $G$ must be a tree. Lemma \ref{lem:treetree}
applies. Let us assume $v(G)>n$. Let $A$ be a set of $n+1$ vertices
spanning a connected subgraph in $G$.  This subgraph must be a
tree. Spoiler plays in $G$ staying all time within
$A$. Lemma~\ref{lem:treetree} applies.

\case3 $G$ is disconnected and has no cycle of length up to $n+1$.

We can assume that every component $H$ of $G$ is a tree for otherwise
Spoiler plays the game on $(T,H)$ staying in $H$, using the strategy
described above.

Suppose first that $G$ has a tree component $H$ such that $H\not\cong
T$ and $v(H)\ge n$.  If $v(H)=n$, let $T'=H$. Otherwise let $T'$ be a
subtree of $H$ on $n+1$ vertices. Spoiler plays the game on $(T,T')$
staying in $T'$ and applying the strategy of Lemma \ref{lem:treetree}
(with $T$ and $T'$ interchanged and perhaps with $n+1$ in place
of~$n$).

Suppose next that all components of $G$ are trees of order less
than~$n$. In the first move Spoiler selects a median $x$ of $T$. Let
Duplicator respond with a vertex $x'$ in a component $T'$ of $G$. If
in the sequel Duplicator makes a move outside $T'$, then Spoiler wins
by Lemma~\ref{lm:path}. As long as Duplicator stays in $T'$, Spoiler
follows the strategy of Lemma \ref{lem:treetree}.

Finally, it remains to consider the case that $G$ has a component $T'$
isomorphic to~$T$. Spoiler plays in $G$. In the first move he
selects a vertex $x'$ outside $T'$. Let $x$ denote Duplicator's
response in $T$.  Starting from the second move Spoiler plays the game
on $(T,T')$ according to Lemma \ref{lem:treetree}, where $x$ is
considered colored in a color absent in $T'$.

Our description of Spoiler's strategy is complete.\epf

It is not clear what the asymptotics of $\DT_0(n,l)$ is. We could not
even rule out the possibility that $\DT_0(n,l)=(\frac12+o(1))\,
\frac{l\ln n}{\ln l}$. 
 \comment{Also, it would be interesting to know $\DT_i$
for other small $i$, such as $i=1$ or $i=2$.}%

The similar method shows that $\DT_0(n,l)=\Theta(\ln n)$ if $l\ge 2$
is constant and $\DT_0(n,l)=\Theta(l)$ if $\frac{\ln n}{\ln l}=O(1)$
but the exact asymptotics seems difficult to compute. 

Using these results, one can show that the upper bounds in
Theorems~\ref{th:RandomTree} and~\ref{th:giant} apply to $D_1(G)$,
that is, there are strategies for Spoiler requiring at most one
alternation. It is not clear whether 0 alternations is possible
here. One of a few places that seem to require an alternation is
establishing that $\phi$ is a bijection: Spoiler may be forced to
start in one of the graphs, while later (for example, when showing that
$T_x\not\cong T'_{x'}$) he may need to swap graphs.

\end{document}

\bibliography{oleg,general,misc,graph,ex,random,enum}
\end{document}